\documentclass[12pt]{article}

\usepackage{amsmath}
\usepackage{amssymb}
\usepackage{amsthm}
\usepackage{bbm}
\usepackage{amscd}
\usepackage{mathabx}
\usepackage[margin=1in]{geometry}
\DeclareFontFamily{U}{wncy}{}
\DeclareFontShape{U}{wncy}{m}{n}{<->wncyr10}{}
\DeclareSymbolFont{mcy}{U}{wncy}{m}{n}
\DeclareMathSymbol{\Sh}{\mathord}{mcy}{"58} 
\newcommand{\calpJ}{\widehat{\mathcal{J}}_{\mathbb{T}_q}}
\newcommand{\pJ} {\widehat{J}_\mathbb{T}}
\newcommand{\calpJQ}{\widehat{\mathcal{J}_{\mathbb{T}_Q}}}
\newcommand{\CC}{\check{H}}

\newtheorem{theorem}{Theorem}[section]
\newtheorem{corollary}[theorem]{Corollary}
\newtheorem{lemma}[theorem]{Lemma}

\title{Freeness of $p$-adically Completed Modular Jacobians over a Hecke Algebra}
\author{John Yu}

\begin{document} 

\maketitle

\section{Introduction}
\subsection{Motivation}
One of the great achievements of modern number theory is Wiles' proof of Fermat's Last Theorem, which crucially relies on an $R = \mathbb{T}$ theorem identifying a certain Hecke ring as a universal deformation ring. The technical underpinning of this isomorphism are Taylor-Wiles systems, which have since been axiomatized and generalized by Diamond and Fujiwara \cite{F} \cite{D}. The original application for the Taylor-Wiles system introduced an auxiliary set of primes at which additional deformation is allowed and a close examination of the ring theoretic properties of Hecke algebras. Taylor-Wiles systems have since been generalized and demonstrated to be a robust framework for tackling a wide range of problems in number theory. In this paper, we set up a Taylor-Wiles system consisting of appropriately chosen $p$-adically completed modular Jacobians over their corresponding Hecke algebras. 
 
\subsection{The main theorem}
Let $K$ be a finite extension of $\mathbb{Q}_p$ with integer ring $\mathcal{O}_K$. Fix a prime $p > 2$, pick a prime $q \equiv 1 \mod p$, and choose a positive integer $N$ such that $(pq,N) = 1$. We fix the prime $p$ throughout this paper, but later $q$ will vary as a Taylor-Wiles prime. We use $\mathcal{O}$ (which is independent of $\mathcal{O}_K$) throughout this paper to denote a discrete valuation ring finite flat over $\mathbb{Z}_p$. Denote by $X'$ and $Y$ the compactified modular curves over $\mathbb{Q}$ associated to the $\Gamma_1(Nq)$ and $\Gamma_1(N) \cap \Gamma_0(q)$ moduli problem of elliptic curves respectively, so that $\text{Gal}(X'/Y) = (\mathbb{Z}/q\mathbb{Z})^\times$. Let $X/Y$ be the intermediate curve with Galois group $\Delta_q = \langle \delta \rangle$, the $p$-Sylow quotient of $(\mathbb{Z}/q\mathbb{Z})^\times$ with a fixed generator $\delta$.  We write $\Gamma_q$ and $\Gamma$ for the congruence subgroup corresponding to $X$ and $Y$ respectively, so that $\Gamma \backslash \Gamma_q$ is canonically isomorphic to $\Delta_q$. 

Denote by $\mathcal{J}$ and $J$ Jacobians over $\mathbb{Q}$ of $X$ and $Y$ respectively and write $\widehat{\mathcal{J}}$ and $\widehat{J}$ for their $p$-adic completions. We fix the notation $X$ and $Y$ to refer to these modular Jacobian. Here, the $p$-adic completion of an abelian variety $A$ over $\mathbb{Q}_\ell$ given by the sheafification of the presheaf $K \mapsto A(K) \otimes_\mathbb{Z} \mathcal{O}$ for $K$ an fppf extension of $\mathbb{Q}_\ell$ (see [H15, \S2] for details of the sheafification process). In this paper, we restrict to finite extensions $K$ of $\mathbb{Q}_p$, in which case this definition coincides with the projective limit $\hat{A}(K) = \varprojlim (A(K)/p^nA(K)) \otimes_\mathbb{Z} \mathcal{O}$ [H15, Lemma 2.1]. From this point on until the appendix, we impose this condition so that $K$ is a $p$-adic field of finite degree over $\mathbb{Q}_p$.

We later also require the notion of $p$-adically completed sheaves evaluated on the \'etale site over a finite field extension $k$ of $\mathbb{F}_p$. Suppose generally that $A$ is an abelian variety over a number field $L$ with good reduction at the prime $p$ and denote by $\mathcal{O}_L$ to be the integer ring of $L$ with residue field $k$. We then let $A_{/\mathcal{O}_L}$ denote the Neron model of $A$ over $\mathcal{O}_L$ and write $\widetilde{A}$ for the special fiber of $A_{/\mathcal{O}_L}$ at $p$. Then, for an \'etale extension $k'$ of the finite field $\mathbb{F}_p$, we define $\hat{A}(k') = \varinjlim \widetilde{A}(\kappa) \otimes \mathcal{O}$ for $\kappa$ running over all finite etale sub-extensions of $k$ inside $k'$. 

We can give an explicit description following \cite[pg 228]{H} of our $p$-adically completed sheaves evaluated on the \'etale site of $\mathbb{Q}_p$. For a finite extension $K$ of $\mathbb{Q}_p$ and an abelian variety $A/\mathbb{Q}_p$, we have the decomposition 

\begin{equation} \label{eq:1.2.1}
A(K) \cong \mathbb{Z}_p^{[K:\mathbb{Q}_p]\cdot\text{dim}(A)} \oplus A[p^\infty](K) \oplus A^{(p)}(K)
\end{equation}

\noindent where $A^{(p)}$ denotes the prime-to-$p$ torsion submodule of $A$ and $A[p^\infty](K)$ is the $p$-torsion subgroup of $A(K)$. The evaluation $\hat{A}(K)$ yields 

\begin{equation} \label{eq:1.2.2}
(A/A^{(p)} \otimes_\mathbb{Z} \mathcal{O})(K) = \mathcal{O}^{[K:\mathbb{Q}_p]\cdot\text{dim}(A)} \oplus (A[p^\infty] \otimes \mathcal{O})(K)
\end{equation}

This is naturally a $\mathcal{O}$-module, and although there is not necessarily an $\mathcal{O}_K$-module structure, $\hat{A}(K)$ at least contains an $\mathcal{O}_K \otimes_\mathbb{Z} \mathcal{O}$-submodule of finite index, which is the formal exponential image of the Lie algebra of $A$ tensored over $\mathbb{Z}$ with $\mathcal{O}$.

We consider these objects as modules over the $p$-adic Hecke algebras $\mathbbm{h} \subset \text{End}_{\mathcal{O}}(\widehat{J})$ and $\mathbbm{h}_q \subset \text{End}_{\mathcal{O}}(\mathcal{\widehat{J}})$ generated by the Hecke operators over $\mathcal{O}$. Fix a non-Eisenstein local ring $\mathbb{T} \subset \mathbbm{h}$ corresponding to the Galois representation $\rho$ valued in $GL_2(\mathbb{T})$ such that $\text{trace}(\rho(\text{Frob}_\ell)) = T_\ell$ for primes $\ell \nmid Np$ (by non-Eisenstein, we mean a local ring that gives rise to an irreducible residual Galois representation). Let $\mathbb{F}$ be the residue field $\mathbb{T}/\mathfrak{m}_\mathbb{T} \mathbb{T}$. We use $\bar{\rho}$ to denote the residual representation valued in $GL_2(\mathbb{F})$ obtained by the reduction mod $\mathfrak{m}_\mathbb{T}$, the maximal ideal of $\mathbb{T}$. By residual irreducibility, this representation is uniquely determined as $\text{trace}(\bar{\rho}(\text{Frob}_\ell))$ is the image of the Hecke operator $T_\ell$ in $\mathbb{T}$ for rational primes $\ell \nmid Np$.  We later make a choice of $q$ as $q$-distinguished primes specified in a Taylor-Wiles system and by their choice of $q$, we can choose a local ring of $\mathbbm{h}_q$ that projects isomorphically onto $\mathbb{T}$, which we denote $\mathbb{T}_q$. The goal of this article is to then give a description of the $\mathbb{T}_q$-module structure of the $p$-completed Mordell-Weil group $\widehat{\mathcal{J}}_{\mathbb{T}_q}(K) = (\widehat{\mathcal{J}} \otimes_{\mathbbm{h}_q} \mathbb{T}_q)(K)$ cut out by $\mathbb{T}_q$ by examining the $\Delta_q$-action on $\calpJ(K)$, given by the diamond operators. 

Before stating our first result, we first recall some technical conditions from Galois deformation theory that we will impose. Let $D_\ell$ denote the decomposition subgroup at the prime $\ell$. We say our local Galois representation $\bar{\rho}: \text{Gal}(\bar{\mathbb{Q}}/\mathbb{Q}) \to GL_2(\mathbb{F})$ is \emph{$p$-distinguished} if $\bar{\rho} |_{D_p}$ has the form $\bar{\rho} |_{D_p} \sim \begin{pmatrix} \bar{\chi_1} & \ast \\ 0 & \bar{\chi_2} \end{pmatrix}$, for characters $\bar{\chi_1} \neq \bar{\chi_2}$. The representation $\bar{\rho}$ is $q$-distinguished if $\bar{\rho}$ is unramified at $q$ and $\bar{\rho} |_{D_q} \sim \begin{pmatrix} \alpha_q & 0 \\ 0 & \beta_q \end{pmatrix}$ for $\alpha_q \neq \beta_q$. In other words, $\bar{\rho}(\text{Frob}_q)$ has two distinct eigenvalues. We say that $\bar{\rho}$ is \emph{$p$-ordinary} if $\bar{\rho} |_{D_p} \sim \begin{pmatrix}\bar{\chi_1} & \ast \\ 0 & \bar{\chi_2} \end{pmatrix}$ for $\bar{\chi_2}$ unramified. We say that $\bar{\rho}$ is \emph{flat} at $p$ if there exists a finite flat group scheme $\mathcal{G}$ over $\mathbb{Z}_p$ with $\mathbb{T}/\mathfrak{m}_{\mathbb{T}}\mathbb{T} \subset \text{End}(\mathcal{G})$ and $\bar{\rho} |_{D_p}$ is realized on $\mathcal{G}(\bar{\mathbb{Q}}_p) \cong (\mathbb{T}/\mathfrak{m}_\mathbb{T}\mathbb{T})^2$ as a Galois module. 

We will later impose additional restrictions on our residual representation $\bar{\rho}$ as part of a Galois deformation problem key in setting up the Taylor-Wiles system in our application. These are some of the necessary conditions that guarantee the existence of a deformation ring $R$, due essentially to Mazur \cite{M}. These technical conditions will be recalled in section 4. 

Finally, we assume throughout this paper that the residual representation $\bar{\rho}$ is absolutely irreducible when restricted to $\mathbb{Q}(\sqrt{(-1)^{(p-1)/2}p})$. This assumption is necessary to guarantee the existence of the Taylor-Wiles primes, which we recall later following Wiles.

The first result needed to prove the main theorem of this paper relates the two $p$-adically completed Jacobians $\pJ$ and $\calpJ$ via the $\Delta_q$-action: 

\begin{theorem} [Theorem A]
\label{thmA} Assume that $\mathbb{T}$ is non-Eisenstein and $K$ is a $p$-adic field of finite degree over $\mathbb{Q}_p$. Suppose that the local representation $\bar{\rho}$ is $q$-distinguished. Assume also that $\pJ(K)$ is torsion-free. We have an isomorphism $\widehat{J}_\mathbb{T}(K)^\ast \cong \calpJ(K)^\ast/(\delta-1)\calpJ(K)^\ast$ where $\calpJ(K)^\ast$ and $\pJ(K)^\ast$ denote the $\mathcal{O}$ dual of $\calpJ(K)$ and $\pJ(K)$ respectively.
\end{theorem}

We later describe explicit sufficient conditions for the torsion-freeness of $\pJ(K)$. 

One can think about $\pJ(K)$ as a fixed base object sitting under a moving family of the objects $\calpJ(K)$, which vary as one adds primes $q \equiv 1 \mod p$ to the level. Under certain assumptions, we show that the $\mathcal{O}[\Delta_q]$-module structure of $\widehat{\mathcal{J}}_\mathbb{T}$ is amenable to an application of Taylor-Wiles systems. In particular, we show that $\calpJ(K)$ is $\mathcal{O}[\Delta_q]$-free and by careful choice of the primes $q$ added to the level of $\calpJ$ following the original argument by Taylor and Wiles, we arrive at the main theorem of the paper:

\begin{theorem} [Theorem B]
\label{thmB}Let the assumptions be as in Theorem A. Suppose additionally that either (1) $\bar{\rho}$ is ordinary and $p$-distinguished or (2) flat. Suppose also that primes $\ell \mid N$ satisfy one of Wiles' conditions (A), (B), or (C) \cite{W} (see section 4.3 below for details). Under these assumptions, $\pJ(K)$ is free of finite positive rank $[K:\mathbb{Q}_p]$ over $\mathbb{T}$.
\end{theorem}

We remark that the conditions imposed at $p$ in the hypothesis of the theorem are equivalent to the Serre weight $k$ of $\bar{\rho}$ satisfying $2 \le k \le p-1$ \cite[pg. 215]{DFG}. Additionally, the hypothesis that $\pJ(K)$ is torsion-free may be replaced with the condition that the $p$-adic cyclotomic character does not arise as a quotient of $\bar{\rho}|_{D_p}$.

The conditions (A), (B), and (C) are recalled in section 4 and are necessary to ensure the deformation theory required for the Taylor-Wiles system works. As we shall see, we define an error term $e_K$ controlling the relationship between $\pJ$ and $\calpJ$. We show in the appendix that the analagous error term vanishes in other situations over $\ell$-adic fields and number fields, which may be of use apart from our application to the main theorem.

\subsection{Organization}
Let us give an overview of this paper's organization.

Section 2 examines the relationship between the two abelian \'etale sheaves $\pJ$ and $\calpJ^{\Delta_q}$. Recall that $\calpJ(K) = \hat{\mathcal{J}} \otimes_{\mathbbm{h}} \mathbb{T}_q$ and $\Delta_q$ is the $p$-Sylow quotient of $(\mathbb{Z}/q\mathbb{Z})^\times$, so that $\calpJ^{\Delta_q}$ is the subsheaf fixed by the diamond operators. In particular, we obtain a abelian \'etale sheaf isomorphism $\pJ \cong \calpJ^{\Delta_q}$, where $\calpJ^{\Delta_q}$ is the subsheaf fixed by the action of the diamond operators. 

This is achieved by examining the Hecke correspondence action on the level of Picard schemes and translating the correspondence action to an explicit one on \v{C}ech cohomology. We show that the kernel and cokernel of a natural map $\widehat{J} \to \widehat{\mathcal{J}}^{\Delta_q}$ can be embedded into \v{C}ech cohomology groups, on which the induced action of the Hecke operators $T_\ell$ act as multiplication by $\ell+1$, implying that the kernel and cokernel of the aforementioned map is Eisenstein. By cutting out our $p$-adically completed modular Jacobians with a non-Eiesenstein local ring, we arrive at the desired isomorphism.

As a consequence, we obtain an exact sequence of abelian \'etale sheaves $$0 \to \pJ \to \calpJ \to \alpha(\calpJ) \to 0$$ where $\alpha = \delta - 1$. 

In section 3, we turn to the question of the exactness of the above exact sequence when we evaluate at a local field $K$. In particular, we define an error term $e_K$, whose vanishing measures the exactness of the above exact sequence evaluated at $K$. 

The main point of this section is to prove that key sequence of \'etale sheaves evaluated at $K$ is exact. This is important in providing an exact control result relating the $p$-adically completed Mordell-Weil groups $\pJ(K)$ and $\calpJ(K)$.  

In section 4, we set up a Taylor-Wiles system $(\mathbb{T}_Q, \widehat{\mathcal{J}}_{\mathbb{T}_Q}(K))_Q$ where $Q$ ranges an infinite set of finite sets of primes $q \equiv 1 \mod p$ following the original work of Wiles. 

We begin by proving the $\mathcal{O}[\Delta_q]$-freeness of $\calpJ(K)^\ast$ where $^\ast$ denotes the $\mathcal{O}$-dual. A crucial ingredient needed for this is the vanishing of the error terms $e_K$ defined in section 2. When the $e_K$ vanish, we obtain an exact sequence $$0 \to \alpha(\calpJ)(K)^\ast \to \calpJ(K)^\ast \to \pJ(K)^\ast \to 0$$ which with some commutative algebra arguments gives the desired freeness result. 

We also briefly review the Taylor-Wiles apparatus and explain how it applies to our setup, concluding with a proof of Theorem B.

We conclude in section 5 with some remarks involving the global situation and questions for future exploration.

\emph{Acknowledgements. }The author would like to thank Haruzo Hida for suggesting this problem and his endless guidance throughout this project. 

\section{The Eisenstein Kernel}

Let $\pi: X \to Y$ be the projection of modular curves as in \S1.2. We have an induced natural map on the $p$-adically completed modular Jacobians $i: \widehat{J} \to \widehat{\mathcal{J}}$, which we consider as a map of abelian \'etale sheaves, and an induced action of $\Delta_q$ on $\calpJ$. A natural question is to ask how $\widehat{J}$ and $\widehat{\mathcal{J}}^{\Delta_q}$ are related. In this section, we show that the natural map between them has Eisenstein kernel and cokernel, so upon tensoring by a non-Eisenstein local component $\mathbb{T}$ of the $p$-adic Hecke algebra, we obtain the following:

\begin{theorem}\label{eis}If $\mathbb{T}$ is non-Eisenstein, we have an isomorphism $\pJ \cong \calpJ^{\Delta_q}$ as sheaves in the \'etale site over $\mathbb{Q}_p$.
\end{theorem}
We prove this in the next subsection following some general results regarding \v{C}ech cohomology. The main point is to interpret the Hecke action on the kernel and cokernel in question in terms of \v{C}ech cohomology groups. We compute explicitly that the action of $T_\ell$ on \v{C}ech cohomology is multiplication by $\ell+1$ and that the kernel and cokernels naturally embed into certain \v{C}ech cohomology groups. 

\subsection{\v{C}ech cohomology}
We begin by interpreting the Hecke correspondence action on Picard groups using \v{C}ech cohomology and demonstrating the relevant technical results. For this section, we work in a general setting. Most of the technical results in this section are found in \cite{H}. We restate some results here as our applications and the context differ. In \cite{H}, the primary application is control of $U_p$-operators, whereas the application here is to explicitly describe the correspondence action of $T_\ell$ operators on cohomology for primes $\ell$ different than $p$. 

Let $\pi: \mathcal{X} \to \mathcal{Y}$ denote a finite flat degree $d$ covering of geometrically reduced, proper varieties over a base field $k$. Let $U \subseteq \mathcal{Y} \times_k \mathcal{Y}$ be a correspondence on $\mathcal{Y}$ with projections $\pi_1$ and $\pi_2$ and denote $U_\mathcal{X}$ the pullback correspondence on $\mathcal{X} \times_k \mathcal{X}$ with projection maps $\pi_{1,\mathcal{X}}$ and $\pi_{2,\mathcal{X}}$. 

For some $\alpha \in H^0(\mathcal{X}, \mathcal{O}_\mathcal{X}^\times)$, consider $\alpha_U := \pi_{1,\mathcal{X}}^\ast \alpha \in H^0(U_\mathcal{X}, \mathcal{O}_{U_\mathcal{X}})$. Multiplication by $\pi_{1,\mathcal{X}}^\ast\alpha$ has a degree $d$ characteristic polynomial $P(T)$ and define the norm $N_U(\alpha_U)$ to be $P(0)$. We set the $U$ action on $H^0(\mathcal{X},\mathcal{O}_\mathcal{X}^\times)$ to be $U(\alpha) = N_U(\alpha_U)$.

To give the correspondence action of $U$ on \v{C}ech cohomology, we next describe the induced action of $U$ on $H^0(\mathcal{X}, \mathcal{O}_{\mathcal{X}^{(s)}}^\times)$. Consider the iterated product $\pi_{i,\mathcal{X}^{(s)}} = \pi_{i,\mathcal{X}} \times_\mathcal{Y} \cdots \times_\mathcal{Y} \pi_{i,\mathcal{X}} : U_\mathcal{X}^{(s)} \to \mathcal{X}^{(s)}$ for $i=1,2$, where $U_\mathcal{X}^{(s)}$ the $s$-fold product of $U_\mathcal{X}$ over $\mathcal{Y}$. Denote by $U_j = U_\mathcal{X}^{(j)} \times_\mathcal{Y} \mathcal{X}^{(s-j)}$ and consider the two projections $\pi_{1,j}^{(s)}$ and $\pi_{2,j}^{(s)}$ to $U_\mathcal{X}^{(j-1)} \times_\mathcal{Y} \mathcal{X}^{(s-j+1)}$. These induce a correspondence $U_j \subseteq (U_\mathcal{X}^{(j-1)} \times_\mathcal{Y} \mathcal{X}^{(s-j+1)}) \times_\mathcal{Y} (U_\mathcal{X}^{(j-1)} \times_\mathcal{Y} \mathcal{X}^{(s-j+1)})$. The norm map $N_{U_s} = N_{\pi_{2,\mathcal{X}^{(s)}}}: \pi_{2, \mathcal{X}^{(s)},\times} \mathcal{O}_{U_s}^\times \to \mathcal{O}_{\mathcal{X}^{(s)}}^\times$ factors through the norm maps as $N_{U_s} = N_s \circ N_{s-1} \circ \cdots \circ N_1$, where $N_j$ is the norm map with respect to $U_j \to U_{j-1}$.

In sufficiently nice situations, we obtain an explicit action of the Hecke operators acting on \v{C}ech cohomology: 

\begin{lemma} $\pi: \mathcal{X} \to \mathcal{Y}$ is a finite flat morphism of geometrically reduced, proper schemes over a field $\text{Spec}(k)$. Suppose that $\mathcal{X}$ and $U_\mathcal{X}$ are proper schemes over $k$ such that

1. $U_\mathcal{X}$ is geometrically reduced, and for each geometrically connected component $\mathcal{X}^\circ$ of $\mathcal{X}$, its pullback to $U_\mathcal{X}$ via $\pi_{2,\mathcal{X}}$ is connected

2. $(f \circ \pi_{2,\mathcal{X}})_\ast \mathcal{O}_{U_\mathcal{X}} = f_\ast \mathcal{O}_\mathcal{X}$ 

We have that 

1. $\pi_2: U \to \mathcal{Y}$ has constant degree $d = \text{deg}(\pi_2)$, the action of $U$ on $H^0(\mathcal{X}, \mathcal{O_\mathcal{X}}^\times)$ is multiplication by $d$.

2. The action of $U^{(s)}$ on $\CC^0(\mathcal{X}, \mathcal{O}_{\mathcal{X}^{(s)}}^\times)$ is multiplication by $d$

3. The action of $U$ on $\CC^s(\mathcal{X}/\mathcal{Y}, H_0(\mathbb{G}_{m,\mathcal{Y}}))$ is multiplication by $d$. 

\end{lemma} 

\emph{Proof}: Since $\pi$ is proper, we have by the assumptions that $H^0(U_\mathcal{X}, \mathcal{O}_\mathcal{X}) = H^0(\mathcal{X},\mathcal{O}_\mathcal{X}) = k^{\left| \pi_0(\mathcal{X}) \right|}$. For a section $\alpha_U \in H^0(U_\mathcal{X}, \mathcal{O}_{U_\mathcal{X}}) = H^0(\mathcal{X},\mathcal{O}_\mathcal{X})$ and so $N_U(\alpha_U) = \alpha_U^{\text{deg}(\pi_2)}$, which gives the first assertion. For the second, we have that $(f \circ \rho_1)_\ast \mathcal{O}_{U_1} = (f \circ \pi_{2,\mathcal{X}})_\ast(\mathcal{O}_\mathcal{X}) \otimes_{\mathcal{O}_\mathcal{Y}} f_\ast\mathcal{O}_\mathcal{X} \otimes_{\mathcal{O}_\mathcal{Y}} \cdots \otimes_{\mathcal{O}_\mathcal{Y}} f_\ast \mathcal{O}_\mathcal{X}$. Therefore, the hypotheses above are satisfied for $U_1$ and we deduce the second assertion. The third follows immediately from the second. \qed

\subsection{Application to $p$-adically completed Jacobians}
The \v{C}ech to derived spectral sequence yields the following commutative diagram with exact rows:
$$\begin{CD}
\CC^1(X/Y, H^0(Y, \mathbb{G}_{m, Y})) @>>> \text{Pic}_{Y/K} @>>> \CC^0(X/Y, \text{Pic}_{Y/K}) @>>> \CC^2(X/Y, H^0(Y, \mathbb{G}_{m,T})) \\
                               @AAA     @AAA     @AAA     @AAA \\
                             C' @>>> J @>>> \CC^0(X/Y, \mathcal{J}) @>>> C
\end{CD}$$
We observe first that since $X/Y$ is Galois that $\CC^0(X/Y, \mathcal{J})$ is just group cohomology $H^0(\Delta_q, \mathcal{J}) = \mathcal{J}^{\Delta_q}$. Since $T_\ell$ acts on the \v{C}ech groups $\CC^n(X/Y, \mathcal{O}_X^\times)$ via multiplication by $\ell+1$, we conclude that after completing $p$-adically (i.e. tensoring by $\mathcal{O}$ over $\mathbb{Z}$) and tensoring by a non-Eisenstein local component of the Hecke algebra $\mathbb{T}$ the desired sheaf isomorphism $\pJ \cong \calpJ^{\Delta_q}$. 

\section{Vanishing of the error term}
We regard as above $\pJ$ and $\calpJ$ as \'etale shaves over the \'etale site of $\mathbb{Q}_p$ and let $K$ be a finite extension of $\mathbb{Q}_p$. In this way, evaluating $\pJ$ and $\calpJ$ gives the $p$-adically completed Mordell-Weil group $\varprojlim J_{\mathbb{T}}(K)/p^{n}J_{\mathbb{T}}(K) \otimes_\mathbb{Z} \mathcal{O}$ and $\varprojlim \mathcal{J}_{\mathbb{T}_q}(K)/p^n\mathcal{J}_{\mathbb{T}_q}(K) \otimes_\mathbb{Z} \mathcal{O}$ respectively. As a result of the previous section, we obtain an exact sequence of abelian \'etale sheaves for a non-Eisenstein local ring $\mathbb{T}$ 

$$0 \to \pJ \to \calpJ \xrightarrow[]{\alpha} \calpJ \xrightarrow[]{\nu} \pJ \to 0$$ 

\noindent or when examining the first three arrows, we consider the short exact sequence 

$$0 \to \pJ \to \calpJ \to \alpha(\calpJ) \to 0$$ 

\noindent where $\alpha = \delta - 1$ and $\nu = \sum_{\sigma \in \Delta_q} \sigma$. 

For $K$ a finite extension of $\mathbb{Q}_p$, define 

$$e_K:=\alpha(\calpJ)(K)/\alpha(\calpJ(K))$$ 

\noindent as the abelian group quotient, which measures the exactness of the corresponding module sequence evaluated at $K$. The goal of this section is to prove that the error term $e_K$ vanishes for any choice of $K$ finite over $\mathbb{Q}_p$.

Before proceeding, we first show that it is enough to show that $e_K$ vanishes for $\mathcal{O}=\mathbb{Z}_p$. Suppose that  $\bar{\rho}$  has values in  $GL_2(\mathbb{F})$ for a finite extension $\mathbb{F}$ of $\mathbb{F}_p$ and let $\mathcal{O}$  be the unramified extension of  $\mathbb{Z}_p$  with residue field  $\mathbb{F}$, i.e.  the Witt vector ring with coefficients in  $\mathbb{F}$. Let $L$ be a number field with integer ring $\mathcal{O}_L$  such that $\mathcal{O}_L \otimes_\mathbb{Z} \mathbb{Z}_p \cong \mathcal{O}$.  For $A_{\mathcal{O}_L}:=A\otimes_\mathbb{Z} \mathcal{O}_L$, where $A$ is either $J$ or $\mathcal{J}$, the sheaf defined by $K \mapsto A_{\mathcal{O}_L}(K) \otimes_\mathbb{Z} \mathbb{Z}_p$  is isomorphic to  $\hat{A}$, where $\hat{A}$ is the $p$-adic completion defined in the introduction. Note that  $A_{\mathcal{O}_L} \cong A^{[\mathbb{F}:\mathbb{Q}]}$  as an abelian variety.  Moreover, for a complete local ring $S$ finite over $\mathbb{Z}_p$  with residue field $\mathbb{F}$, $S$  is automatically an $\mathcal{O}$-algebra by the universality of Witt ring, and hence  $S \otimes_{\mathbb{Z}_p} \mathcal{O}=S^m$  for  $m=[\mathbb{F}:\mathbb{F}_p]$  (as $S \otimes \mathcal{O}/\mathfrak{m}_S = \mathbb{F} \otimes_{F_p} \mathbb{F} =\mathbb{F}^m$).  Thus, if we can prove control for  $\mathbb{T}$  (and $\mathbb{T}_q$)-component of $\hat{J}$  and $\hat{\mathcal{J}}$ for  $\mathbb{Z}_p$, it is valid for the unramified $\mathcal{O}$.  For more general $\mathcal{O}$, we can further extend scalars from the Witt vector ring to $\mathcal{O}$.  Therefore, in the rest of this section, we assume $\mathcal{O}=\mathbb{Z}_p$.

\subsection{Vanishing at $p$}

We discuss first the $p$-adically completed \'etale sheaf of an abelian variety $A$ over a finite field of characteristic $p$. Let $A$ be an abelian variety over the finite field $\mathbb{F}_p$ and write, for an \'etale extension $k$ of $\mathbb{F}_p$, $k = \varinjlim_n F_n$ for finite etale subextensions of $F_n$ of $k$. For an etale extensions $k$ on the \'etale site of $\mathbb{F}_p$, we define $\displaystyle \hat{A}(k) = \varinjlim_{n} (A(F_n) \otimes \mathbb{Z}_p) = \bigcup_{n} A(F_n) \otimes \mathbb{Z}_p$. For $k$ a finite extension of $\mathbb{F}_p$, we simply have $\hat{A}(k) = A(k) \otimes \mathbb{Z}_p$. The sheafification of this presheaf yields the $p$-adic completion of $A$. As tensoring with $\mathbb{Z}_p$ over $\mathbb{Z}$ is exact, this procedure transforms exact sequences of abelian varieties over $\mathbb{F}_p$ to exact sequences of \'etale sheaves on the \'etale site over $\mathbb{F}_p$.

\begin{lemma} 
\label{lang} For $A$ an abelian variety over a finite field $k$ of characteristic $p$, consider the $p$-adically completed Mordell-Weil group $\widehat{A}(k)$. Under these hypotheses, we have the vanishing of the group cohomology group $H^1(G_k,\widehat{A}(\overline{k})) = 0$, where $G_k = \text{Gal}(\overline{k}/k)$, the absolute Galois group of $k$.
\end{lemma}

\emph{Proof}: Write $q=p^n$ where $n$ is the degree of $k$ over $\mathbb{F}_p$. From the usual Lang's theorem, we have that the morphism $\sigma: A(\bar{k}) \to A(\bar{k})$ where $\sigma(x)=x^{-1}\cdot x^q$ is a surjection. The kernel of $\sigma$ is precisely $A(k)$ and looking at the corresponding long exact sequence of cohomology and using the fact that tensoring with $\mathbb{Z}_p$ commutes with formation of Galois cohomology (as it is flat over $\mathbb{Z}$), we conclude that $H^1(G_k,A(\overline{k}))=0$ and and by taking the $p$-primary part (or by tensoring with $\mathbb{Z}_p$) that $H^1(G_k,\widehat{A}(\overline{k})) = 0$. \qed

\begin{lemma}
\label{vanishing} For $K$ of residual characteristic $p$, $e_K$ vanishes, i.e. $\alpha(\pJ(K))=\alpha(\pJ)(K)$.
\end{lemma}

\emph{Proof}: In the following for an abelian variety $A$ over $\mathbb{Q}$, the notation $\widehat{A}(\mathbb{F})$ is taken to mean evaluating the $p$-adic completion of the special fiber of the Neron model of $A$ at the $\mathbb{F}$ in the manner described above (in applications, $A$ will be either $\mathcal{J}$ or $J$). Consider the commutative diagram with exact columns, where the kernel of reduction $A^\circ(K)$ is the formal group at the identity. 
$$\begin{CD}
0 @>>> \pJ^\circ(K) @>>> \calpJ^\circ(K) @>>> \alpha(\calpJ)^\circ(K) @>>> 0 \\
  @.                @VVV                   @VVV                           @VVV @.\\
0 @>>> \pJ(K) @>>> \calpJ(K) @>>> \alpha(\calpJ)(K) @>>> 0 \\
@. @VVV @VVV @VVV @. \\
0 @>>> \pJ(\mathbb{F}) @>>> \calpJ(\mathbb{F}) @>>> \alpha(\calpJ)(\mathbb{F}) @>>> 0 
\end{CD}$$
\noindent As we work on the \'etale site over $\mathbb{Q}_p$, we obtain a description of $\pJ(K)$ and $\calpJ(K)$ explicitly as $\mathbb{Z}_p$-modules and the mod $p$ reductions and apply Lemma \ref{lang} to deduce that $H^1(\text{Gal}(\bar{\mathbb{F}}/\mathbb{F}), \pJ(\bar{\mathbb{F}})) = 0$. From the long exact sequence arising from the bottom row, we deduce exactness of the bottom sequence. As $\pJ$ and $\calpJ$ have good reduction at $p$, the kernel of reduction (i.e. the formal group at the identity) is formally smooth and the sequence 

$$0 \to \pJ^\circ \to \calpJ^\circ \to \alpha(\calpJ)^\circ \to 0$$

\noindent is exact as an etale sheaf. By formal smoothness, write $\pJ^\circ(K)$ and $\alpha(\calpJ)^\circ(K)$ as $\text{Spf}(\mathbb{Z}_p[[x_1,\dots,x_i]])$ and $\text{Spf}(\mathbb{Z}_p[[y_1,\dots,y_j]])$. The sheaf exactness tells us that the middle term of the top exact sequence is given by the formal spectrum of the completed tensor product $\mathbb{Z}_p[[x_1,\dots,x_i]] \widehat{\otimes}_{\mathbb{Z}_p} \mathbb{Z}_p[[y_1,\dots,y_j]]$. We thus conclude that the top row sequence of $p$-adically completed Mordell-Weil groups is also exact. Finally, the nine lemma yields exactness of the middle row, as desired. \qed

\begin{lemma}
\label{controlseq} 
When $K$ is a finite extension of $\mathbb{Q}_p$, the following sequence is exact:
$$0 \to \alpha(\calpJ)(K) \to \calpJ(K) \to \nu(\calpJ(K)) \to 0$$
\end{lemma} 

\emph{Proof: }Immediate from the vanishing of $e_K$ and left exactness of sheaf evaluation. \qed

\section{The Taylor-Wiles System}
The goal of this section to prove Theorem B. We begin by proving the $\mathcal{O}[\Delta_q]$-freeness of $\calpJ(K)$ and subsequently generalizing this freeness to $p$-adically completed modular Jacobians $\calpJQ(K)$, where $\calpJQ$ will range over an infinite set of finite sets $Q$ to be defined precisely later. As in the original Taylor-Wiles paper, these freeness results are crucial in setting up the Taylor-Wiles system we use. We also briefly review Taylor-Wiles systems and deformation theory of Galois representations as needed for our applications. We conclude the section with some remarks on the torsion-free assumption of our main theorems.

\subsection{Freeness over $\mathcal{O}[\Delta_q]$}
We wish to prove the $\mathcal{O}[\Delta_q]$-freeness of the $\calpJ(K)$. To do so, we recall a key result due to Taylor and Wiles \cite[Theorem 2]{TW}, demonstrating the freeness of Hecke algebras $\mathbb{T}_q$ over the group ring $\mathcal{O}[\Delta_q]$. 

\begin{lemma}
\label{heckefree} $\mathbb{T}_q$ is $\mathcal{O}$-free over $\mathbb{T}$ of rank equal to the $\mathcal{O}$ rank of $\mathbb{T}$.
\end{lemma}

\emph{Proof}. This is the content of Theorem 2 of \cite{TW}. \qed

For the following result, we specialize briefly to $\mathcal{O} = \mathbb{Z}_p$ and describe how to obtain the analogous statement for general $\mathcal{O}$ after.

\begin{lemma}
\label{jacfree} Suppose $\pJ(K)$ is $\mathbb{Z}_p$-free and $\text{dim}_{\mathbb{T}/\mathfrak{p} \otimes \mathbb{Q}_p} \widehat{J}_\mathbb{T}(K)[\mathfrak{p}] \otimes \mathbb{Q}_p = r$ for all $\mathfrak{p} \in \text{Spec}(\mathbb{T})(\overline{\mathbb{Q}_p})$. Then, we have the bound $\text{dim}_{\mathbb{T}_q/\mathfrak{Q} \otimes \mathbb{Q}_p} \calpJ(K)[\mathfrak{Q}] \otimes \mathbb{Q}_p \le r$ for all $\mathfrak{Q} \in \text{Spec}(\mathbb{T}_q)(\overline{\mathbb{Q}_p})$. If we assume further that $\text{dim}_{\mathbb{T}/\mathfrak{Q} \otimes \mathbb{Q}_p} \calpJ(K)[\mathfrak{Q}] \otimes \mathbb{Q}_p = r$ for all $\mathfrak{Q} \in \text{Spec}(\mathbb{T}_q)(\overline{\mathbb{Q}_p})$, then $\calpJ(K)$ is a free $\mathbb{Z}_p[\Delta_q]$-module of rank $rR$, where $R$ is the $\mathbb{Z}_p$-rank of $\mathbb{T}$.
\end{lemma}

\emph{Proof}. Since $\pJ(K)$ is $\mathbb{Z}_p$ torsion-free, we have that $H^0(\text{Gal}(\overline{\mathbb{Q}_p}/K), \bar{\rho})=0$. Since $\calpJ(K)$ shares the same residual representation $\bar{\rho}$ as $\pJ(K)$, we conclude that $\calpJ(K)$ is also torsion free and thus $\mathbb{Z}_p$-free by examining the decomposition (1.2.1) for $\calpJ(K)$. 

Recall the sheaf exact sequence

$$0 \to \pJ \to \calpJ \xrightarrow[]{\alpha} \calpJ \xrightarrow[]{\nu} \pJ \to 0$$ 

\noindent where upon evaluating at our $p$-adic field $K$. We split this 4 term sequence into the following exact sequences:

$$0 \to \pJ(K) \to \calpJ(K) \to \alpha(\calpJ(K)) \to 0$$
$$0 \to \alpha(\calpJ(K)) \to \calpJ(K) \to \nu(\calpJ(K)) \to 0$$

\noindent in which every term is $\mathbb{Z}_p$-free. Combining these two exact sequences with an application of Lemma 3.3 and 3.4, we obtain the isomorphism $\pJ(K)^\ast \cong \calpJ(K)^\ast/(\delta - 1)\calpJ(K)^\ast$. 

Since $\pJ(K)$ has $\mathbb{Z}_p$ rank equal to $r \cdot R$, we obtain a surjection $\mathbb{Z}_p[\Delta_q]^{rR} \to \calpJ(K)^\ast$ as $\mathbb{Z}_p[\Delta_q]$-modules. Let $\mathfrak{Q}_0$ be a prime ideal of $\mathbb{Z}_p[\Delta_q]$ with residual characteristic 0. Since $\mathbb{T}_q$ is a $\mathbb{Z}_p[\Delta_q]$-free module of finite rank $R$ (as $\mathbb{T}_q/(\delta - 1)\mathbb{T}_q \cong \mathbb{T}$), we have the equalities

$$R\text{dim}_{\mathbb{T}_q/\mathfrak{Q}_0 \otimes \mathbb{Q}_p} \calpJ(K)^\ast/\mathfrak{Q}_0 \otimes \mathbb{Q}_p \le rR \cdot \text{rank}_{\mathbb{Z}_p}\mathbb{Z}_p[\Delta_q]/\mathfrak{Q}_0 = r\cdot\text{rank}_{\mathbb{Z}_p} \mathbb{T}_q/\mathfrak{Q}_0 = rR$$

On the other hand, we have 

$$\text{dim}_\mathbb{F}(\calpJ(K)^\ast/\mathfrak{Q}_0\calpJ(K)^\ast) \otimes_{\mathbb{Z}_p[\Delta_q]/\mathfrak{Q}_0} \mathbb{F} = \text{dim}_\mathbb{F} \calpJ(K)^\ast \otimes_{\mathbb{Z}_p[\Delta_q]} \mathbb{F} = rR$$

\noindent since $\mathfrak{Q}_0 \subset \mathfrak{m}_{\mathbb{Z}_p[\Delta_q]}$. We conclude that the minimal number of generators of $\calpJ(K)^\ast/\mathfrak{Q}_0\calpJ(K)^\ast$ over $\mathbb{Z}_p[\Delta_q]$ is $rR$ by Nakayama's lemma.

By assumption, we have that $\text{dim}_{\mathbb{T}/\mathfrak{Q} \otimes \mathbb{Q}_p} \calpJ(K)[\mathfrak{Q}] \otimes \mathbb{Q}_p = r$ and so $\calpJ(K)[\mathfrak{Q}] \otimes \mathbb{Q}_p \cong k(\mathfrak{Q})^r$ for the residue field $k(\mathfrak{Q})$. This implies $(\calpJ(K)^\ast/\mathfrak{Q}\calpJ(K)^\ast) \otimes \mathbb{Q}_p \cong k(\mathfrak{Q})^r$ since $k(\mathfrak{Q})$ is self-dual over $\mathbb{Q}_p$ by the trace pairing. We conclude therefore that 

$$\text{dim}_{\mathbb{T}/\mathfrak{Q} \otimes \mathbb{Q}_p}(\calpJ(K)^\ast/\mathfrak{Q}\calpJ(K)^\ast) \otimes \mathbb{Q}_p = r$$

\noindent for any prime $\mathfrak{Q} \in \text{Spec}(\mathbb{T})$ with residual characteristic zero. We also have that

$$(\calpJ(K)^\ast/\mathfrak{Q}_0\calpJ(K)^\ast) \otimes \mathbb{Q}_p \cong \prod_{\mathfrak{Q} \mid \mathfrak{Q}_0} (\calpJ(K)^\ast/\mathfrak{Q}\calpJ(K)^\ast) \otimes \mathbb{Q}_p$$

\noindent for any prime $\mathfrak{Q}_0 \in \text{Spec}(\mathbb{Z}_p[\Delta_q])$ with residual characteristic zero and deduce from this relation that

$$(\calpJ(K)^\ast/\mathfrak{Q_0}\calpJ(K)^\ast) \otimes \mathbb{Q}_p \cong (\mathbb{Z}_p[\Delta_q]/\mathfrak{Q}_0) \otimes_{\mathbb{Z}_p} \mathbb{Q}_p^{rR} = k(\mathfrak{Q}_0)^{rR}$$

\noindent for the residue field $k(\mathfrak{Q}_0 = \text{Frac}(\mathbb{Z}_p[\Delta_q]/\mathfrak{Q}_0)$. Since $\mathbb{Z}_p[\Delta_q]/\mathfrak{Q}$ is a discrete valuation ring, we have that $\calpJ(K)^\ast/\mathfrak{Q}_0\calpJ(K)^\ast \cong (\mathbb{Z}_p[\Delta_q]/\mathfrak{Q}_0)^{rR} \oplus T$ for some finite torsion module $\mathbb{Z}_p[\Delta_q]/\mathfrak{Q}_0$-module $T$. Since $\calpJ(K)^\ast/\mathfrak{Q}_0\calpJ(K)^\ast$ has minimal number of generators equal to $rR$ over $\mathbb{Z}_p[\Delta_q]/\mathfrak{Q}_0$, we must have $T = 0$, and so we obtain 

$$\calpJ(K)^\ast/\mathfrak{Q}_0\calpJ(K)^\ast \cong (\mathbb{Z}_p[\Delta_q])^{rR}$$

\noindent for all $\mathfrak{Q}_0$. This isomorphism is induced by $\pi: \mathbb{Z}_p[\Delta_q]^{rR} \to \calpJ(K)^\ast$ tensored with $\mathbb{Z}_p[\Delta_q]/\mathfrak{Q}_0$ and so $\text{ker}(\pi) \subset \cap_{\mathfrak{Q}_0} \mathfrak{Q}_0^{rR} = 0$, with $\mathfrak{Q}_0$ running over all primes with residual characteristic zero. We conclude therefore that $\calpJ(K)^\ast \cong \mathbb{Z}_p[\Delta_q]^{rR}$ and so $\calpJ(K)^\ast$ is $\mathbb{Z}_p[\Delta_q]$-free of rank $R$. To obtain the freeness of $\calpJ(K)$, we take the $\mathbb{Z}_p$ dual noting that $\mathbb{Z}_p[\Delta_q]$ is a local complete intersection (as we can exhibit an explicit isomorphism $\mathbb{Z}_p[\Delta_q] \cong \mathbb{Z}_p[S]/((1+S)^{p^{e}})$, where $e$ is the $p$-order of the cyclic group $\Delta_q$). \qed

Later, we will show that $r=1$ in our situation. We obtain as a corollary the analagous freeness result over $\mathcal{O}[\Delta_q]$. 

\begin{lemma}
\label{ojacfreeness}Under the same assumptions as Lemma \ref{jacfree}, $\calpJ(K)$ is free of finite rank over $\mathcal{O}[\Delta_q]$.
\end{lemma} 
\emph{Proof. } Writing  $\mathbb{T}_{q/\mathcal{O}}$  and  $\mathbb{T}_{/\mathcal{O}}$ for the local ring of the Hecke algebra over  $\mathcal{O}$, when  $\mathcal{O}=W$, where $W$ is the Witt vector ring with residue field $\mathbb{F}$), we have that $\mathbb{T}_{/\mathcal{O}} \cong \mathbb{T}_{/\mathbb{Z}_p}^m$ and $\mathbb{T}_{q/\mathcal{O}} \cong \mathbb{T}_{q/\mathbb{Z}_p}^m$ for $m=[\mathbb{F}:\mathbb{F}_p]$. As explained at the beginning of Section 3,  $\hat {\mathcal{J}}_{\mathbb{T}_q/\mathcal{O}}=\widehat{\mathcal{J}}_{\mathbb{T}_q/\mathbb{Z}_p}^m$, and hence  we can replace $\mathbb{Z}_p$ in the above argument by  $\mathcal{O}=W$.  For general $\mathcal{O}$ with reisdue field $\mathbb{F}$, one obtain the result by scalar extension from $W$ as $\mathbb{T}_{q/W} \otimes_W \mathcal{O} = \mathbb{T}_{q/\mathcal{O}}$. This proves that the $\mathcal{O}[\Delta_q]$-module $\hat{\mathcal{J}}_{\mathbb{T}_q}(K)$  is  $O[\Delta_q]$-free. \qed

We next let $Q=\{q_1,\dots,q_r\}$ be a finite set of primes $q_i \equiv 1 \mod p$ for which $\bar{\rho}$ is $q_i$-distinguished and unramified at each $q_i$. Let $Z'$ be the compact modular curve at level $\Gamma_1(N \cdot \prod_{q \in Q} q)$ and let $Z$ be the curve intermediate to $Z'$ and $X$ whose Galois group is precisely $\Delta_Q$, the $p$-Sylow quotient of $\text{Gal}(Z'/X) = \prod_{q \in Q} (\mathbb{Z}/q\mathbb{Z})^\times$. Let $\mathbbm{h}_Q$ denote the $\mathcal{O}$-algebra of Hecke operators acting on the $p$-adically completed Jacobian $\widehat{J}$ and let $\mathbb{T}_Q$ denote a local ring of $\mathbbm{h}_Q$ which projects down to $\mathbb{T}$ according to Taylor and Wiles \cite[pg. 555]{TW}. We use $\widehat{\mathcal{J}}_{\mathbb{T}_Q}$ to denote the tensor product $\widehat{\mathcal{J}} \otimes_{\mathbbm{h}_Q} \mathbb{T}_Q$. 

By inducting, we obtain the corollary:

\begin{corollary} 
\label{delQfree} Suppose $\pJ(K)$ is $\mathcal{O}$-free and $\text{dim}_{\mathbb{T}/\mathfrak{p} \otimes \mathbb{Q}_p} \pJ(K)[\mathfrak{p}] \otimes \mathbb{Q}_p = r$ for all $\mathfrak{p} \in \text{Spec}(\mathbb{T})(\overline{\mathbb{Q}_p})$. Then, we have that $\text{dim}_{\mathbb{T}_Q \otimes \mathbb{Q}_p} \widehat{\mathcal{J}}_{\mathbb{T}_Q}(K)[\mathfrak{Q}] \otimes \mathbb{Q}_p \le r$ for all $\mathfrak{Q} \in \text{Spec}(\mathbb{T}_Q)(\overline{\mathbb{Q}_p})$. Assume further that $\text{dim}_{\mathbb{T}_Q \otimes \mathbb{Q}_p} \widehat{\mathcal{J}}_{\mathbb{T}_Q}(K)[\mathfrak{Q}] \otimes \mathbb{Q}_p = r$ for all $\mathfrak{Q} \in \text{Spec}(\mathbb{T}_Q)(\overline{\mathbb{Q}_p})$. We have that $\widehat{\mathcal{J}_{\mathbb{T}_Q}}(K)$ is a free $\mathcal{O}[\Delta_Q]$-module of rank $rR$, where $R$ is the $\mathcal{O}$ rank of $\mathbb{T}$.
\end{corollary}
\emph{Proof}. We induct on the cardinality of the set $Q$. The base case where $\left|Q\right| = 1$ is Lemma \ref{ojacfreeness}. Suppose that the result is true for any $Q$ with $\left|Q\right| = r-1$ and suppose inductively that $Q = \{q_1,\dots,q_r\}$. Write $Q' = Q - \{q_r\}$ and let $\Gamma_Q$ and $\Gamma_{Q'}$ denote the corresponding levels groups, so that $\Gamma_Q/\Gamma_{Q'} \cong \Delta_{q_r} = \langle \delta \rangle$. Using similar arguments from the previous sections switching the roles of $\pJ(K)$ and $\calpJ(K)$ with $\widehat{\mathcal{J}_{\mathbb{T}_Q'}}$ and $\widehat{\mathcal{J}_{\mathbb{T}_Q}}$, we deduce the relation

$$\widehat{\mathcal{J}_{\mathbb{T}_Q}}(K)^\ast/(\delta-1)\widehat{\mathcal{J}_{\mathbb{T}_Q}}(K)^\ast \cong \widehat{\mathcal{J}_{\mathbb{T}_{Q'}}}(K)^\ast$$

\noindent By our inductive step, the right term, $\widehat{\mathcal{J}}_{\mathbb{T}_{Q'}}(K)^\ast$ is $\mathcal{O}[\Delta_{Q'}]$-free of rank $r \cdot R$. As such, we obtain a surjective morphism $\mathcal{O}[\Delta_Q]^{rR} \twoheadrightarrow \widehat{\mathcal{J}_{\mathbb{T}_Q}}(K)^\ast$. As $\mathbb{T}_Q$ has rank $R$ over $\mathcal{O}$, we then follow the argument of Lemma \ref{jacfree} replacing $\pJ(K)$ with $\widehat{\mathcal{J}_{\mathbb{T}_Q'}}$ and $\calpJ(K)$ with $\widehat{\mathcal{J}_{\mathbb{T}_Q}}(K)$. \qed

\subsection{Taylor-Wiles Systems}
We  recall the Taylor-Wiles formalism, following Fujiwara \cite{F}, Diamond \cite{D}, and Hida \cite{HB}:

As consistent with the rest of this paper, $K$ denotes a finite extension of $\mathbb{Q}_p$, $\mathcal{O}$ denotes its integer ring, and $\mathbb{F}$ its residue field. The symbol $\mathfrak{m}_\mathcal{O}$ denotes the maximal ideal of $\mathcal{O}$. 

Let $X$ be a set of finite subsets of primes, which also contains $\emptyset$. We take a pair $(R,M)$, where $R$ is a complete noetherian local $\mathcal{O}$-algebra with the residue field $\mathbb{F}$ and $M$ a finitely generated $R$-module. A Taylor-Wiles system $\{R_Q, M_Q\}_{Q \in X}$ for $(R,M)$ consists of the following data:

\begin{itemize}
\item [(TW1)] For $Q \in X$ and $q \in Q$, $q \equiv 1 \mod p$. We denote by $\Delta_q$ the $p$-Sylow subgroup of $(\mathbb{Z}/q\mathbb{Z})^\times$ and $\Delta_Q$ is defined as $\prod_{q \in Q} \Delta_v$ for $Q \in X$. 

\item [(TW2)] For $Q \in X$, $R_Q$ is a complete noetherian local $\mathcal{O}[\Delta_Q]$-algebras with residue field $\mathbb{F}$ and $M_Q$ is an $R_Q$-module. For $Q = \emptyset$, $(R_\emptyset, M_\emptyset) = (R,M)$. 

\item [(TW3)] A surjection $R_Q/I_Q \to R_Q$ of local $\mathcal{O}$-algebras for each $Q \in X$. Here, $I_Q \subset \mathcal{O}[\Delta_Q]$ denotes the augmentation ideal of $\mathcal{O}[\Delta_Q]$. For $Q = \emptyset$, it is the identity of $R$. 

\item [(TW4)] The homomorphism $R_Q/I_QR_Q \to \text{End}_{\mathcal{O}_\lambda}(M_Q/I_QM_Q)$ factors through $R$ and $M_Q/I_QM_Q$ is isomorphic to $M$ as an $R$-module

\item [(TW5)] $M_Q$ is free of finite rank $\alpha$ as an $\mathcal{O}[\Delta_Q]$-module for a fixed $\alpha$.

\end{itemize}
The main theorem for Taylor-Wiles systems is the following:

\begin{theorem}
\label{twsystem} For a Taylor-Wiles system $\{R_Q,M_Q\}_{Q \in D}$ for $(R,M)$, assume the following conditions.

\begin{itemize}
\item [(1)] \emph{For any integer $n \in \mathbb{N}$, the subset $X_n$ of $X$ defined by $X_n = \{Q \in X \mid q \in Q \implies q \equiv 1 \mod p^n\}$ is an infinite set.}

\item [(2)] \emph{$r = \left|Q\right|$ is independent of $Q \in X$ if $Q \neq \emptyset$}

\item [(3)] \emph{$R_Q$ is generated by at most $r$ elements as a complete local $\mathcal{O}$-algebra for all $Q \in X$.}
\end{itemize}

Under these conditions, (1) $R$ is $\mathcal{O}$ flat and is a local complete intersection ring and (2) $M$ is a free $R$-module. In particular, $R$ is isomorphic to the image $T$ in $\text{End}_{\mathcal{O}} M$. 
\end{theorem}

\subsection{Galois Deformation Theory}

In this section, we recall some of the theory of Galois deformations and formulate the particular deformation problem for our application. Let $\Sigma$ be a finite set of primes including the prime $p$ and write $\mathbb{Q}^\Sigma$ to denote the maximal extension of $\mathbb{Q}$ unramified outside $\Sigma$ and $\infty$. We fix an embedding $\mathbb{Q}$ and $\mathbb{Q}^\Sigma$ in the complex numbers and we fix a choice of decomposition groups $D_\ell$ for all rational primes $\ell$. Let $k$ be a finite field of characteristic $p$ and let  $\bar{\rho_0}: \text{Gal}(\mathbb{Q}^\Sigma/\mathbb{Q}) \to GL_2(k)$ be an irreducible representation. We assume that $\text{det}(\bar{\rho_0})$ is odd. We say that $\rho_0$ is a \emph{deformation} of $\bar{\rho_0}$ if $\rho_0$ is a representation of $\text{Gal}(\mathbb{Q}^\Sigma/\mathbb{Q})$ valued in an Artinian local ring $A$ with residue field $k$ such that the residual representation is equal to $\bar{\rho_0}$. Following Mazur, we consider deformations equivalent up to a conjugation action of elements of the kernel of the reduction map $GL_2(A) \to GL_2(k)$. 

In applications, it is useful to consider deformations of $\bar{\rho}$ that satisfy a range of conditions naturally arising from geometry. A \emph{deformation condition} means a full subcategory of the deformations to $\mathbb{T}$ satisfying axioms (1)-(3) of \cite[pg. 289]{M}. It is a theorem of Mazur utilizing Schlessinger's criterion that a formal deformation condition is relatively representable and that if the residual representation $\bar{\rho}$ is absolutely irreducible, then it is represented by a quotient ring of the larger deformation problem. In our situation, we consider deformations of our residual deformation $\bar{\rho}$ with coefficients in $\mathbb{F}$ from the non-Eisenstein maximal ideal $\mathfrak{m}$, so our deformation problems seek lifts $\rho$ with coefficients in $\mathbb{T}$. We consider the following deformation types:

\begin{itemize}
\item
\begin{enumerate}
\item \emph{Selmer deformations}. We assume that $\bar{\rho}$ is ordinary and that the deformation has a representative $\rho: \text{Gal}(\mathbb{Q}^\Sigma/\mathbb{Q}) \to GL_2(\mathbb{T})$ such that $\rho |_{D_p} \sim \begin{pmatrix} \chi_1 & \ast \\ 0 & \chi_2 \end{pmatrix}$ with $\chi_2$ unramified and $\chi_1 \equiv \chi_2 \mod \mathfrak{m}_\mathbb{T}$ with $\text{det} |_{I_p} = \epsilon\omega^{-1}\chi_1\chi_2$ for the cyclotomic character $\epsilon$ and $\omega$ is of order prime to $p$ satisfying $\omega \equiv \epsilon \mod p$. 

\item \emph{ordinary deformations}. We impose the same conditions as in (1) but with no determinant conditions

\item \emph{strict deformations}. The same as in (2), but imposed when $\bar{\rho}$ is not semisimple and not flat. In this case, we demand that $\bar{\chi_1}\bar{\chi_2}^{-1} = \omega$. 

\end{enumerate}

Remark: As we work with weight 2 modular forms, cases (2) and (3) do not arise in our applications.

\item
 \emph{Flat at $p$ deformations}. We assume in this case that each deformation $\rho$ to $GL_2(\mathbb{T})$ has the property that for any quotient $\mathbb{T}/\mathfrak{\alpha}$ of finite order $\rho |_{D_p} \mod \mathfrak{\alpha}$ is the Galois representation associated to a finite flat group scheme over $\mathbb{Z}_p$. 

\end{itemize}

Additionally, the following types of minimal deformation conditions will be imposed at the primes of bad reduction that are necessary for the Taylor-Wiles method to work. These conditions were originally imposed by Wiles. There are three behaviors of the local Galois deformation we distinguish at primes $\ell$ where the representation $\bar{\rho}$ is ramified. These are not necessarily mutually exclusive conditions.

\begin{itemize}

\item[(A)] $\bar{\rho}|_{D_\ell} = \begin{pmatrix} \alpha & \ast \\ 0 & \beta \end{pmatrix}$ for a suitable choice of basis with $\alpha, \beta$ unramified characters satisfying $\alpha\beta^{-1} = \omega$ and the fixed space of $I_\ell$ of dimension 1.

\item[(B)] $\bar{\rho} |_{I_\ell} = \begin{pmatrix} \alpha & 0 \\ 0 & 1 \end{pmatrix}$, with $\alpha \neq 1$ for a suitable choice of basis

\item[(C)] $H^1(\mathbb{Q}_\ell, W) = 0$ where $W$ is the adjoint representation of $\bar{\rho}$, where the adjoint representation $W$ is the space of matrices $\text{Mat}_{2\times 2}(\mathbb{F})$ whose Galois module action is given by $x \mapsto gxg^{-1}$ for $g \in D_\ell$ and $x \in \text{Mat}_{2 \times 2}(\mathbb{F})$. 

\end{itemize}

We impose for primes $\ell$ satisfying (A), (B), or (C) the condition that their deformations to $\mathbb{T}$ satisfy the following conditions respectively:

\begin{itemize}

\item[(A$'$)] $\rho |_{D_\ell} = \begin{pmatrix} \psi_1 & \ast \\ 0 & \psi_2 \end{pmatrix}$ with $\psi_1, \psi_2$ unramified and $\psi_1\psi_2^{-1} = \epsilon$

\item[(B$'$)] $\rho |_{I_\ell} = \begin{pmatrix} \chi_\ell & 0 \\ 0 & 1 \end{pmatrix}$ for a suitable choice of basis, with $\chi_\ell$ of prime order to $p$

\item[(C$'$)] $\text{det}(\rho |_{I_\ell})$ is of order prime to $p$

\end{itemize}
Following Wiles notation, $\mathcal{M}$ will denote the ramified primes satisfying (A), (B), or (C). 

We now describe the deformation problem that will be of concern to us. Let $Q$ be a set of primes $q \equiv 1 \mod p$ such that $\bar{\rho}$ is unramified and $q$-distinguished at each $q \in Q$. We let $\mathcal{J}_Q$ denote the modular Jacobian whose Galois group over $J$ is precisely $\Delta_Q = \prod_{q \in Q} \Delta_q$ and write $\mathbb{T}_Q$ to denote the $p$-adic Hecke algebra of $\mathcal{J}_Q$ projecting to $\mathbb{T}$. In particular, these Hecke algebras give rise to Galois representations with the same residual representation $\bar{\rho}$.

\noindent We let $\mathcal{M}$ be the primes dividing $N$ and let $\Sigma$ denote the union of the primes $\mathcal{M} \cup \{p\} \cup Q$, where we allow $Q = \emptyset$ (the \emph{minimal} case). Let $G_Q$ denote the Galois group $\text{Gal}(\mathbb{Q}^{\Sigma}/\mathbb{Q})$. Let $\text{CNL}_{\mathcal{O}}$ denote the category of complete, Noetherian local rings with the same residue field $\mathbb{F}$ as $\mathcal{O}$, where the morphisms are given by local $\mathcal{O}$-algebra local homomorphisms. We define the functor 

$$\mathcal{D}_{(\cdot,Q,\Sigma,\mathcal{M})}: \text{CNL}_{\mathcal{O}} \to \text{Sets}$$ 

\noindent where $\mathcal{D}_{(\cdot,Q,\Sigma,\mathcal{M})}(B)$ is the set of deformations $\rho: G_Q \to GL_2(B)$ of $\bar{\rho}$ such that

\begin{enumerate}
\item the representation $\rho$ satisfies $\cdot$, where $\cdot$ is either the condition of being ordinary and $p$-distinguished or (2) flat

\item At primes $\ell \in \mathcal{M}$ satisfying (A), (B), (C), the representation $\rho$ is required to satisfy conditions (A$'$), (B$'$), (C$'$) respectively  

\end{enumerate}

Work of Mazur and Ramakrishna gives the existence of a unique deformation ring $R_{(\cdot,Q,\mathcal{M})}$, which is complete, local, and Noetherian, and represents this functor. In our application, the set $\mathcal{M}$ always remains the same as the places where the residual representation $\bar{\rho}$ attached to $\mathbb{T}$ is ramified, i.e. those primes dividing $N$. Additionally, the set $\Sigma$ depends only on $Q$. We therefore drop the $\Sigma$ and $\mathcal{M}$ with this understanding and abbreviate our deformation problem $\mathcal{D}_{(\cdot,Q,\Sigma)}$ to $\mathcal{D}_Q$. We similarly drop the $\cdot$ from the notation with the understanding that we are either in the (1) ordinary and $p$-distinguished or (2) flat case. We abbreviate the representing deformation ring in this case to be $R_Q$. 

Let $R$ be the deformation ring representing the deformation problem $\mathcal{D}_\emptyset$ and let $M$ be the Betti cohomology group $H^1(X(\Gamma_Q), \mathcal{O})$. Work of Wiles shows that there is an infinite family $X$ of finite sets $Q_j$ of primes satisfying the Taylor-Wiles conditions (TW1-5) \cite[pg. 523]{W}. For these sets $Q_j$, Wiles shows that $R_Q \cong \mathbb{T}_Q$. From this point onwards, we utilize the same family of finite sets $X = \{Q_m\}$, but we will in the next section introduce our own Taylor-Wiles system replacing the Betti cohomology groups with the $p$-adically completed Jacobians $\widehat{\mathcal{J}}_{\mathbb{T}_Q}(K)$ that we have been working with throughout this paper. 

\subsection{Freeness of $\pJ(K)$}
\emph{Proof of Theorem B. }Since the residual representation $\bar{\rho}$ is absolutely irreducible when restricted to $\mathbb{Q}(\sqrt{(-1)^{p(p-1)/2}})$, we follow Wiles in \cite[pg. 522]{W} and produce an infinite set $X$ of finite sets $Q_m$ consisting of primes $q \equiv 1 \mod p^m$ and $\left|Q_m\right| = r$ such that each $q \in Q_m$ is $q$-distinguished. Given these $Q_m$, we recall the deformation functor defined in the preceding section

$$\mathcal{D}_{Q_m}: \text{CNL}_{\mathcal{O}} \to \text{Sets}$$

\noindent which is represented by the deformation ring $R_{Q_m}$ and is isomorphic by Wiles to the Hecke algebra $\mathbb{T}_{Q_m}$. 

Let $t_{\pJ^\circ}$ and $t_{\calpJQ^\circ}$ denote the tangent space of modular Jacobians $\pJ$ and $\calpJQ$. By Tate \cite[pg. 169]{T}, there exists a local isomorphism 

$$\text{exp}: t_{\pJ^\circ} \to \pJ^\circ(K)$$

\noindent to the $K$-rational points of the formal group of $\pJ$, which is $\mathcal{O}$-linear and becomes an isomorphism upon tensoring with $\mathbb{Q}_p$. We use here $\text{exp}$, the inverse of the log map that Tate explitly defines in \cite{T}. Recall $Y$ as in section 1.2 and regard $Y$ to be defined over $\mathcal{O}$. As we have the identifications

$$t_{\pJ^\circ} \cong H^1(Y,\mathcal{O}_Y) \otimes_{\mathbbm{h}} \mathbb{T} \cong \mathbb{T}$$

\noindent as $\mathbb{T}$-modules, we obtain an embedding of $\mathbb{T}$ as an $\mathcal{O}$-lattice in a finite index subgroup of $\pJ(K)$. Here, we have used the fact that $\mathbbm{h}$ is $\mathcal{O}$-dual to $H^0(Y,\Omega_{Y/\mathcal{O}})$ \cite[Theorem 3.17]{HB} and that $H^0(Y,\Omega_{Y/\mathcal{O}})$ is $\mathcal{O}$-dual to $H^1(Y,\mathcal{O}_Y)$ by Grothendieck-Serre duality.  Since tensoring with $\mathbb{Q}_p$ yields an isomorphism, we deduce that $\text{dim}_{\mathbb{T}/\mathfrak{p} \otimes \mathbb{Q}_p} \pJ(K)[\mathfrak{p}] = 1$ for all $\mathfrak{p} \in \text{Spec}(\mathbb{T})(\overline{\mathbb{Q}_p})$. By a similar argument using the exponential map, we have an embedding of the Hecke $\mathbb{T}_Q$ in a finite index subgroup of $\calpJQ(K)$. As $\text{rank}_{\mathcal{O}[\Delta_Q]} \mathbb{T}_Q = R$, we conclude that $\text{dim}_{\mathbb{T}_Q/\mathfrak{Q} \otimes \mathbb{Q}_p} \calpJQ(K)[\mathfrak{Q}] \otimes \mathbb{Q}_p = 1$ as well. Applying Corollary \ref{delQfree}, we conclude that $\calpJQ(K)$ is free of finite rank $R$ over $\mathcal{O}[\Delta_Q]$. 

We now set $R = \mathbb{T}$ and $M = \pJ(K)$. With these notations, $(R_{Q_m}, M_{Q_m})$ is a Taylor-Wiles system for $(R,M)$ in the sense of the preceding section. We write $\Delta_{Q_m} = \prod_{q \in Q_m} \Delta_q$ where the $\Delta_q$ are the diamond operators at the appropriate level. We verify each of the following Taylor-Wiles conditions:

\begin{itemize}
\item [(TW1-3)] are clear either by definition or construction.

\item[(TW4)] is a consequence of Theorem A.

\item[(TW5)] is a consequence of the discussion above
\end{itemize}

By Corollary \ref{delQfree} combined with Theorem \ref{twsystem} applied to  $M=\pJ(K)$, we therefore obtain the $\mathbb{T}$-freeness of $\pJ(K)$. To determine the $\mathbb{T}$ rank of $\pJ(K)$, we observe that $\pJ(K) \cong \mathcal{O}^{[K:\mathbb{Q}_p] \cdot R}$, where $R$ is the $\mathcal{O}$-rank of $\mathbb{T}$. Comparing $\mathcal{O}$-ranks, we conclude that the $\mathbb{T}$-rank of $\pJ(K)$ is $[K:\mathbb{Q}_p]$. This concludes the proof of Theorem B. \qed

\subsection{Criteria for torsion-freeness}
We conclude this section with some remarks on the torsion-free assumption of $\pJ(K)$. In particular, we describe criteria to guarantee that $\pJ(K)$ is torsion-free. We distinguish between two cases, namely the case where $\bar{\rho}$ is ordinary following work of Hida and the non-ordinary case following work of Serre. 

When $\bar{\rho}$ is ordinary, let $\Lambda = \mathbb{Z}_p[[T]]$ be the Iwasawa algebra. The Galois action factors through $\mathbb{Z}_p^\times = \Gamma \times \mu_{p-1} \to \Lambda^\times$. Then, $\zeta \in \mu_{p-1}$ maps to $\zeta^{a_\mathbb{T}}$ for some $0 \le a_\mathbb{T} < p-1$. Let $a(p)$ denote the image of the $U(p)$-operator (which lives originally in $\mathbbm{h}$) in $\mathbb{T}$. We have the following result due to Hida \cite[Lemma 4.1]{H}:

\begin{lemma}[Hida]
\label{ordvanish} Suppose $K$ is finite over $\mathbb{Q}_p$ and let $\mathbb{T}$ be a reduced local ring of the $p$-adic Hecke algebra $\mathbbm{h}$ of $J$. Let $f$ be the residual degree of $K$ and suppose either of the following conditions are satisfied:

\begin{itemize}
\item [(a)] \emph{$a_K(\mathbb{T}) \neq 0$ and $a(p)^f \not\equiv 1 \mod \mathfrak{m}_\mathbb{T}$}

\item [(b)] \emph{$\mathbb{T}$ is Gorenstein and $\bar{\rho}(I_p)$ has a non-trivial unipotent element, where $I_p$ is the inertia group of $\text{Gal}(\bar{\mathbb{Q}_p}/K)$}
\end{itemize}

Then, $\pJ(K)$ is torsion-free.
\end{lemma}

On the other hand, in the non-ordinary case, if $\bar{\rho} |_{I_t}$ decomposes into two characters of level two for the tamely ramified inertia subgroup $I_t$, then the representation $\bar{\rho}$ is irreducible over $\mathbb{Q}_p$ \cite{S}. Therefore, if $K$ is linearly disjoint from the splitting field of $\bar{\rho} |_{\text{Gal}(\overline{\mathbb{Q}_p}/\mathbb{Q}_p)}$, it remains irreducible over $K$. In these two cases, $\pJ(K)$ is $p$-torsion free and so $\pJ(K)$ is torsion free. 

\section{Concluding remarks}
For the Taylor-Wiles system to apply, we required the $\mathbb{Z}_p$ freeness of $\pJ(K)$ and $\calpJ(K)$, which occurs only when the residual characteristic of $K$ is $p$. In other cases (i.e. when the residual characteristic of $K$ is $\neq p$, the $K$-rational points of the $p$-adic completions of the modular Jacobians instead form finite abelian groups. Though in these situations we cannot apply the methodology of the final section, the analysis of the error terms still remain valid and thus the error terms $e_K$ vanish in many situations. When the error terms $e_K$ do not vanish, however, we can then produce non-trivial elements of the Tate-Shafarevich group $\Sh_{\pJ}(F)$ for a global field $F$ when $e_{F_v}$ is non-zero for a finite place $v$ of $F$. We hope to study this relationship and the case where $K$ is a number field in the future. 

As the methods used in this paper did not depend crucially on the fact we were in the elliptic modular case, another question for future investigation is whether or not the same arguments can be used in a Shimura curve setting. 

\section{Appendix: The Error Term for $\ell \neq p$}
We maintain the same notation as before, except that $K$ is of residual characteristic $\ell \neq p$. The aim of this appendix is to show that the error terms $e_K$ defined earlier vanish under suitable assumptions. We still work over the \'etale site of $\mathbb{Q}_\ell$, in which case our notion of $p$-adic completion follows the same definition as before by Hida. 

\subsection{Vanishing in the good reduction case}
Suppose $\ell$ is the residual characteristic of $K$ and that $\ell \nmid Nqp$.

\textbf{Lemma A.1. }\emph{As $\ell \nmid Npq$, $\mathcal{J}$ has good reduction at $\ell$. If the residual characteristic of $K$ is $\ell \nmid Npq$ and the ramification index of $\mathcal{O}_K$ is $e(\ell) < \ell -1$, then $e_K$ vanishes.}

Pf. If the hypotheses are satisfied, we can apply Theorem 4 of BLR Chapter 7.5 Theorem 4 so we obtain an exact sequence of abelian \'etale sheaves over $\mathcal{O}_K$:

\begin{center}
$0 \to \pJ/\mathcal{O}_K \to \calpJ/\mathcal{O}_K \to \alpha(\calpJ)/\mathcal{O}_K \to 0$
\end{center}

By Lang's theorem, $H^1(\mathbb{F}, \pJ) = 0$ for the residue field $\mathbb{F}$ of $\mathcal{O}_K$ and so we obtain an exact sequence

\begin{center}
$0 \to \pJ(\mathbb{F}) \to \calpJ(\mathbb{F}) \to \alpha(\calpJ)(\mathbb{F}) \to 0$
\end{center}

Since $K$ is an $\ell$-adic local field, we have that $\pJ(K) = \pJ[p^\infty](K) = \pJ[p^\infty](\mathbb{F}) = \pJ(\mathbb{F})$ and similarly $\calpJ(K) = \calpJ(\mathbb{F})$ since outside the level, $\pJ$ and $\calpJ$ have good reduction. We conclude therefore that $e_K = 0$.  

\subsection{Vanishing in the bad reduction case}
When the residual characterstic of $K$ divides $N$, we make the assumption that the representation $\bar{\rho}\mid_{D_\ell}$ does not contain the $p$-adic cyclotomic character as a quotient.

\textbf{Lemma A.2. }\emph{Suppose the residual characteristic $\ell$ of $K$ divides $N$ and that the representation $\bar{\rho}\mid_{D_\ell}$ does not contain the mod $p$ $p$-adic cyclotomic character as a quotient. Then, $e_K$ vanishes.}

\emph{Proof}: Since either the cyclotomic character does not appear as a quotient of the local Galois representation, we have that $H_0(\mathbb{Q}_\ell,J[p]_\mathbb{T}(-1)) = 0$. As a result, we have that $H_0(\mathbb{Q}_\ell,T_p\widehat{J}(-1)) = 0$.

Given our exact cohomology sequence

\begin{center}
$0 \to \widehat{J}_{\mathbb{T}}(K) \to \calpJ(K) \to \alpha(\calpJ)(K) \to H^1(\pJ, K) \to \cdots$
\end{center}

By local Tate duality, we have that $H^1(K, \pJ) \cong$ $^t\pJ(K)\;\widecheck{  }$ and similarly $H^1(K, \calpJ) \cong$ $^t\calpJ\widecheck{  }\;$. Denote by $Z = \text{ker}(^t\pJ(K)\;\widecheck{  } \to$ $^t\calpJ\widecheck{  }\;)$. We wish to show that this vanishes, giving the desired exactness. But as $H_0(K, T_pJ_{\mathbb{T}}(-1))$ is dual to $^t\widehat{J_{\mathbb{T}}}(K)\;\widecheck{  }\;$, we obtain the desired vanishing. 

\emph{Remark: }When $K = \mathbb{Q}_\ell$, we can give an explicit descripton in terms of the automorphic representations. Let $\rho = \widehat{\otimes}_\ell \pi_\ell$ be the automorphic representation corresponding to $\overline{\rho}$. Then, the condition that $\bar{\rho}$ does not contain the cyclotomic character as a quotient translates to the following situations \cite{G}, where $\chi_p$ is the $p$-adic cyclotomic character:

(a) $\pi_\ell$ is \emph{principal}, so that $\pi_\ell = \pi(\alpha, \beta)$, and we demand that $\chi_p \mod p \notin \{\alpha \mod p,\beta \mod p\}$

(b) $\pi_\ell$ is \emph{special}, so that $\pi_\ell = \sigma(\eta, \left| \cdot \right|_\ell^{-1} \eta)$ for a character $\eta$ of $\mathbb{Q}_\ell^\times$, and we demand that $\chi_p\mod p\not\in\{|\cdot|_l\eta\mod p,\eta\mod p\}$

(c) $\pi_\ell$ is \emph{supercuspidal} and $\ell \neq 2$ (the $\ell = 2$ case is exceptional for technical reasons, so we do not consider it here), and the corresponding local Galois representation modulo $p$ does not contain $\left| \cdot \right|_\ell$


\begin{thebibliography}{9}
 
\bibitem[BLR90]{BLR} 
Sigfried Bosch, Werner Lutkebohmert, Michel Raynaud. 
\textit{Neron Models}.  
Springer-Verlag (1990).

\bibitem[D97]{D}
Fred Diamond
\textit{The Taylor-Wiles construction and multiplicity one}.
Inventiones mathematicae 128 (1997). 

\bibitem[DFG04]{DFG}
Fred Diamond, Matthias Flach, Li Guo.
\textit{The Tamagawa number conjecture of adjoint motives of modular forms}.
Annales scientifiques de l'Ecole Normale Superieure, Serie 4: Volume 37 (2004) no 5 (2004)

\bibitem[F06]{F}
Kazuhiro Fujiwara 
\textit{Deformation Rings and Hecke Algebras in the Totally Real Case}. 
arxiV (2006).

\bibitem[G75]{G}
Stephen S. Gelbart
\textit{Automorphic Forms on Adele Groups}.
Princeton University Press and University of Tokyo Press (1975). 

\bibitem[H15]{H} 
Haruzo Hida. 
\textit{Limit Mordell-Weil Groups and Their $p$-adic Closure}. 
Documenta Mathematica: Extra Volume Merkurjev (2015).

\bibitem[MFGC]{HB}
Haruzo Hida.
\textit{Modular Forms and Galois Cohomology}.
Cambridge Studies in Advanced Mathematics (2000).

\bibitem[Maz97]{M}
Barry Mazur.
\textit{An introduction to the deformation theory of Galois representations}.
Modular Forms and Fermat's Last Theorem, Springer-Verlag, (1997).

\bibitem[S87]{S}
Jean-Pierre Serre 
\textit{Sur les repr\'esentations modulaires de degr\'e 2 de Gal($\overline{\mathbb{Q}}/\mathbb{Q}$)}. 
Duke Mathematics Journal, Vol. 54, No. 1 (1987).

\bibitem[T67]{T}
John Tate 
\textit{$p$-Divisible Groups}. 
Proceedings of a Conference on Local Field (1967).

\bibitem[TW95]{TW}
Richard Taylor and Andrew Wiles
\textit{Ring-theoretic properties of certain Hecke algebras}.
Annals of Mathematics, 141 (1995).

\bibitem[W95]{W}
Andrew Wiles 
\textit{Modular Elliptic Curves and Fermat's Last Theorem}. 
Annals of Mathematics, 142 (1995).


\end{thebibliography}
\end{document}